\documentclass[12pt]{amsart}
\newcommand{\pof}{\noindent{\em Proof: }}
\usepackage{color}

\newcommand{\s}[1]{\mathcal{#1}}

\newtheorem{Thm}{Theorem}[section]
\newtheorem{Def}[Thm]{Definition}
\newtheorem{Rem}[Thm]{Remark}
\newtheorem{Lem}[Thm]{Lemma}

\newtheorem{Prop}[Thm]{Proposition}

\numberwithin{equation}{section}

 \fontsize{.5cm}{.5cm}\selectfont\sf

\medskip

\begin{document}

\title[Quantum special linear group]
{Dual canonical bases for the quantum special linear group and
invariant subalgebras} \thanks{Supported by the Australian
Research Council and Chinese National Natural Science Foundation
project number:10471070}
\author{Hechun Zhang}
\address{Department of Mathematical Sciences, Tsinghua University, Beijing,
100084, P. R. China} \email{hzhang@math.tsinghua.edu.cn}
\author{R. B. Zhang}
\address{ School of Mathematics and Statistics, University of Sydney,
Australia} \email{rzhang@maths.usyd.edu.au}
\begin{abstract} A string basis is constructed for each
subalgebra of invariants of the function algebra on the quantum
special linear group. By analyzing the string basis for a
particular subalgebra of invariants, we obtain a ``canonical
basis''  for every finite dimensional irreducible
$U_q({\mathfrak{sl}}(n))$-module. It is also shown that the
algebra of functions on any quantum homogeneous space is generated
by quantum minors.
\end{abstract}
\maketitle

\section{Introduction}
The theory of crystal and canonical bases was initiated and
developed by Kashiwara \cite{k1, k} and Lusztig \cite{lu1, lu}.
Many remarkable features are known for canonical bases, among
which the positivity property seems to have the deepest
implications. Variations of canonical bases were also introduced
and investigated in the literature. In \cite{du}, a dual basis of
the canonical basis of the modified quantum enveloping algebra of
type $A$ was investigated under the name global IC basis. In
\cite{le2}, the dual canonical basis of $U_q({\mathfrak n}_+)$ was
constructed by using the so-called quantum shuffles, and was shown
to be related to the representation theories of Hecke algebras and
quantum affine algebras \cite{lnt} (see also \cite{mtz}). In
\cite{bz1} Berenstein and Zelevinsky conjectured a {\em
multiplicative property}, which states that two dual canonical
basis elements $b_1, b_2$ of $U_q({\mathfrak n}_+)$ $q$-commute if
and only if $b_1 b_2= q^m b$ for some $b$ in the dual canonical
basis and some integer $m$. The conjecture was studied by using
Hall algebra techniques in reference \cite{r}. It was observed
that a large portion of the dual canonical basis enjoyed the
multiplicative property, see also \cite{cr}. However,
counter-examples are given to Berenstein and Zelevinsky's
conjecture in \cite{le1} by finding some so-called imaginary
vectors. Reference \cite{bz1} also introduced the notion of string
bases, which was designed to realize the so-called good bases. See
\cite{mathieu} for the existence of good bases.

In \cite{zhanghc} the dual basis of the canonical basis of the
modified quantum enveloping algebra was constructed and the case
of type $A$ was extensively investigated. It was shown that the
dual canonical basis is invariant under the multiplication of
certain $q$ central elements. More importantly, it was shown in
\cite{jz} that the product of any two dual canonical basis
elements is a ${\mathbb Z}_+[q, q^{-1}]$ combination of elements
of the dual canonical basis, establishing multiplicative
positivity of the dual canonical basis (see Theorem
\ref{positivity}).

In this note we use the dual canonical basis constructed in
\cite{zhanghc} for the algebra of functions ${\s O}_q(SL(n))$ on
the quantum special linear group to study its subalgebras of
invariants. Recall that ${\s O}_q(SL(n))$ admits two natural left
actions of the quantized universal enveloping algebra
$U_q({\mathfrak{sl}}(n))$ respectively corresponding to the left
and right translations familiar in the classical context of Lie
groups. Given any subalgebra ${\s C}$ of $U_q({\mathfrak{sl}}(n))$
which is at the same time a two-sided co-ideal, one can show
\cite{gz} that the subspace of ${\s O}_q(SL(n))$ invariant under
left translation with respect to ${\s C}$ forms a subalgebra. Such
invariant subalgebras of ${\s O}_q(SL(n))$ include the algebras of
functions on quantum homogeneous spaces \cite{nmy, gz} as special
cases.

A string basis (in the sense of \cite{bz1}) is constructed for
each subalgebra of invariants in this paper. In the case when ${\s
C}$ is $U_q({\mathfrak n}_-)$, we obtain from the string basis a
``canonical basis" for every irreducible finite dimensional
$U_q({\mathfrak{sl}}(n))$-module (see Section \ref{inv} for
details). When the subalgebra of invariants defines a quantum
homogeneous space, we show that it is generated by quantum minors
of some particular form.

The arrangement of the paper is as follows. In Section 2, we
recall the construction of the dual canonical basis of ${\s
O}_q(SL(n))$ in \cite{zhanghc}.  In Section 3, we introduce a new
type of bases for ${\s O}_q(SL(n))$, which are suitable for
studying the actions of the Kashiwara  operators. Section 4
contains the main results of the paper. We first show that any
subalgebra of invariants is spanned by a subset of the dual
canonical basis of ${\s O}_q(SL(n))$. [This result is stated in a
more precise manner in Theorem \ref{main}.] We then derive from
this fact the results discussed in the last paragraph.

We mention that our approach to the canonical bases of finite
dimensional irreducible $U_q({\mathfrak{sl}}(n))$-modules is very
different from the usual approach taken in canonical basis theory.
Here we first realize every finite dimensional irreducible
$U_q({\mathfrak{sl}}(n))$-module by using the quantum Borel-Weil
theorem, then construct a canonical basis for it by exploiting
Theorem \ref{main}.

\section{The construction of the basis of ${\s O}_q(M(n))$}

Through out the paper, the base field is ${\mathbb Q}(q)$, where
$q$ is an indeterminate over the rational numbers. The coordinate
algebra  ${\s O}_{q}(M(n))$ of the quantum matrix is an
associative algebra, generated by elements
$x_{ij},i,j=1,2,\cdots,n$, subject to the following defining
relations:
\begin{eqnarray}\label{relations}x_{ij}x_{ik}&=&q^2
x_{ik}x_{ij} \text{ if } j<k,\\\nonumber x_{ij}x_{kj}&=&q^2x_{kj}x_{ij}
\text{ if }i<k,\\\nonumber
x_{ij}x_{st}&=&x_{st}x_{ij}\text{ if } i>s, j<t,\\\nonumber
x_{ij}x_{st}&=&x_{st}x_{ij}+(q^2-q^{-2})x_{it}x_{sj} \text{ if }
i<s, j<t.\end{eqnarray} For any matrix $A=(a_{ij})_{i,j=1}^n\in
M_n({\mathbb Z}_+)$ ( ${\mathbb Z}_+=\{0,1,\cdots\}$) we define a
monomial $x^A$  by
\begin{equation} x^A=\Pi_{i,j=1}^nx_{ij}^{a_{ij}},\end{equation}
where the factors are arranged in the lexicographic order on
$I(n)=\{(i,j)\mid i,j=1,\dots,n\}$. More explicitly,
$$x^A=x_{11}^{a_{11}}x_{12}^{a_{12}}\cdots x_{1n}^{a_{1n}}
x_{21}^{a_{21}}\cdots x_{2n}^{a_{2n}} \cdots x_{nn}^{a_{nn}}.$$
Since the algebra ${\s O}_q(M(n))$ is an iterated Ore extension
(see, e.g., pp11-14 in \cite{GW} for a discussion on iterated Ore
extensions), the set $\{x^A|A\in M_n({\mathbb Z}_+\}$ is a basis
of the algebra ${\s O}_q(M(n))$.

Let us recall the construction of the dual canonical basis in
\cite{zhanghc}. From the defining relations (\ref{relations}) of
the algebra ${\s O}_q(M(n))$, it is easy to show the following
lemma.

\begin{Lem}
\begin{enumerate}
\item The mapping
\begin{eqnarray*}^-:x_{ij} \mapsto x_{ij}, &&
q\mapsto q^{-1}\end{eqnarray*} extends to an algebra
anti-automorphism of ${\s O}_{q}(M(n))$ regarded as an algebra
over ${\mathbb Q}$.

\item The mapping
\begin{equation*}\sigma: x_{ij}\mapsto x_{ji}\end{equation*}
extends to an algebra automorphism of ${\s O}_{q}(M(n))$ regarded
as an algebra over $K={\mathbb Q}(q)$.
\end{enumerate}
\end{Lem}

\begin{Rem} For the relation between the bar action in the
present form and the bar action in the papers
by Kashiwara, see \cite{zhanghc}.\end{Rem}

Given $A=(a_{ij})_{n\times n}\in M_n({\mathbb Z}_+)$, we define
the row sum $ro(A)$ and the column sum $co(A)$ of the matrix,
respectively, by
\begin{eqnarray*}
ro(A)&=&(\sum_j a_{1j},\cdots,\sum_ja_{nj})=(r_1,r_2,\cdots,r_n), \\
co(A)&=&(\sum_j a_{j1},\cdots,\sum_j a_{jn})=(c_1,c_2,\cdots,c_n).
\end{eqnarray*}

Let $M=(m_{ij})\in M_{n}({\mathbb Z}_+)$.  If $m_{ij}m_{st}\ge 1$
for two pairs of indices $i,j$ and $s,t$ satisfying $i<s, j<t$, we
define a new matrix $M^\prime =(m_{uv}^\prime)\in M_{n}({\mathbb
Z}_+)$ with \begin{eqnarray*}
m_{ij}^\prime=m_{ij}-1, && m_{st}^\prime=m_{st}-1,\\
m_{it}^\prime=m_{it}+1, && m_{sj}^\prime=m_{sj}+1,\\
m_{uv}^\prime=m_{uv},  && \mbox{ for all other entries}.
\end{eqnarray*}
We say that the matrix $M^\prime$ is obtained from the matrix $M$
by a $2\times 2$ sub-matrix transformation. Using this we may
define a partial order on the set $M_{n}({\mathbb Z}_+)$ such that
$M< N$ if $M$ can be obtained from $N$ by a sequence of $2\times
2$ sub-matrix transformations. Note that we have $ro(M)=ro(N)$,
$co(M)=co(N)$.

From the defining relations (\ref{relations}) of the algebra ${\s
O}_q(M(n))$, we have
\[\overline{x^A}=E(A)x^A+\sum_{B}c_{BA}x^B,\]
where
\[
E(A)=q^{-2(\sum_i\sum_{j>k}a_{ij}a_{ik}+\sum_i\sum_{j>k}a_{ji}a_{ki})}\]
and $c_{BA}\in{\mathbb Z}[q,q^{-1}]$, which is nonzero only if
$B<A$.

Let $D(A)=q^{-\sum_i \sum_{j>k}a_{ij}a_{ik}-\sum_i
\sum_{j>k}a_{ji}a_{ki}}$ and let $x(A)=D(A)x^A.$ Set
$${\s L}^*=\oplus_{A\in M_n({\mathbb Z}_+)}{\mathbb Z}[q]x(A).$$
In \cite{zhanghc}, a basis of the algebra ${\s O}_q(M(n))$ is
constructed which has the following property:

\begin{Thm}\label{basis}
There is a unique basis $B^*=\{b(A)|A\in M_n({\mathbb Z}_+)\}$ of
${\s L}^*$ such that:
\begin{enumerate}
\item $\overline{b(A)}=b(A)$ for all $A$.
\item $b(A)=x(A)+\sum_{B<A} h_{BA}x(B)$ where $h_{BA}\in q{\mathbb Z}$$[q]$.
\end{enumerate}
\end{Thm}

\begin{Def} The basis $B^*$ of the algebra ${\s O}_q(M(n))$  is called
the dual canonical basis
of ${\s O}_q(M(n))$.\end{Def}

The quantum determinant ${\det}_q$ is defined in the usual way by:
\begin{equation} {\det}_q=
\Sigma_{\sigma\in S_n}(-q^2)^{l(\sigma)}x_{1\sigma(1)}x_{2\sigma(2)}
\cdots x_{n\sigma(n)}.\end{equation}
It is known that $det_q$ is a central and
group-like element of the algebra  $ {\s O}_q(M(n))$.

Let $m\le n$ be a positive integer. Let $I=\{i_1,i_2,\cdots,i_m\}$
and $J=\{j_1,j_2,\cdots,j_m\}$ be any two subsets of
$\{1,2,\cdots,n\}$, each having cardinality $m$. By examining the
defining relations (\ref{relations}) one can see that the
subalgebra of $ {\s O}_q(M(n))$ generated by the elements $x_{ij}$
with $i\in I, j\in J$, is isomorphic to ${\s O}_q(M(m))$, so we
can talk about its determinant. Such a determinant is called a
quantum minor, and will be denoted by ${\det}_q(I,J)$.

In \cite{zhanghc}, it was also proved that
\begin{Prop} \label{q-minors} All of the quantum minors are dual canonical
basis elements. Moreover, the basis $B^*$ is invariant under the
multiplication of the quantum determinant. Furthermore,
$\sigma(B^*)=B^*$.\end{Prop}

The algebra of functions ${\s O}_q(SL(n))$ on the quantum special
linear group will be identified with ${\s O}_q(M(n))/\langle
\det_q-1 \rangle$. For convenience of reference, we denote by
\begin{equation}
T: {\s O}_q(M(n)) \rightarrow {\s O}_q(M(n))/\langle \det_q - 1
\rangle \label{T}
\end{equation} the canonical map. Then the image
$\bar{B^*}=T(B^*)$ of the set $B^*$ forms a canonical basis of ${\s
O}_q(SL(n))$, which is called the dual canonical basis of ${\s
O}_q(SL(n))$.  In \cite{jz}, it was proved that the dual canonical
basis $\bar{B^*}$ enjoys a remarkable property, the positivity of
the multiplication.

\begin{Thm} \label{positivity}
For any $b_1, b_2\in \bar{B^*}$, we have
$$b_1b_2=\sum_{b\in \bar{B^*}, n\in{\mathbb Z}}
c_{b_1, b_2, b, n}q^nb,$$
where $c_{b_1, b_2, b,n}\in
{\mathbb Z}_+$ are zero except for finitely many $b, n$.\end{Thm}

The following  formulas are  needed for the computation involving
the Kashiwara operators. Define that
$[m]=\frac{q^{-4m}-1}{q^{-4}-1}$.

\begin{Lem}For any $i<k, j<l$,
\begin{eqnarray*}
x_{kl}x_{ij}^s=x_{ij}^sx_{kl}+(q^{2-4s}-q^{-2})x_{ij}^{s-1}x_{il}x_{kj},\\
x_{kl}^sx_{ij}=x_{ij}x_{kl}^s+(q^{2-4s}-q^2)x_{il}x_{kj}x_{kl}^{s-1}.
\end{eqnarray*} Furthermore,
\begin{equation}\label{power}
(x_{ij} x_{kl}-q x_{il}x_{kj})^s=\sum_{m=0}^s (-q^2)^m
\begin{pmatrix}s\\ m\end{pmatrix}_{q^{-4}} q^{4m(m-s)}x_{ij}^{s-m}x_{il}^m
x_{kj}^mx_{kl}^{s-m}.\end{equation}
\end{Lem}

\pof We prove the third equation only, using induction on $s$. The
case with $s=1$ is clearly true. Assume that the equation holds
for $s$. Then
\begin{eqnarray*}&&(x_{ij}x_{kl}-qx_{il}x_{kj})^{s+1}\\
&=&(x_{ij}x_{kl}-qx_{il}x_{kj})\sum_{m=0}^s
(-q^2)^m\begin{pmatrix}s\\ m\end{pmatrix}_{q^{-4}}
q^{4m(m-s)}x_{ij}^{s-m}x_{il}^m x_{kj}^m x_{kl}^{s-m}.
\end{eqnarray*}
The right hand side of the equation can be express as
\begin{eqnarray*} && \sum_{m=0}^s (-q^2)^m \begin{pmatrix}s\\
m\end{pmatrix}_{q^{-4}}
q^{4m(m-s)}x_{ij}x_{kl}x_{ij}^{s-m}x_{il}^m x_{kj}^mx_{kl}^{s-m}\\
&+&\sum_{m=0}^s (-q^2)^{m+1} \begin{pmatrix}s\\
m\end{pmatrix}_{q^{-4}} q^{4(m+1)(m-s)}x_{ij}^{s-m}x_{il}^{m+1}
x_{kj}^{m+1}x_{kl}^{s-m}\\ &=&\sum_{m=0}^s (-q^2)^{m+1}
\begin{pmatrix}s\\m\end{pmatrix}_{q^{-4}}
q^{4(m+1)(m-s)}x_{ij}^{s-m}x_{il}^{m+1} x_{kj}^{m+1}x_{kl}^{s-m}\\
&+&\sum_{m=0}^s (-q^2)^m \begin{pmatrix}s\\
m\end{pmatrix}_{q^{-4}} q^{4m(m-s)}
x_{ij}[x_{ij}^{s-m}x_{kl}+(q^{-4(s-m)+2}-q^2)
x_{ij}^{s-m-1}x_{il}x_{kj}] \\ \nonumber &\times&x_{il}^m x_{kj}^m
x_{kl}^{s-m}.
\end{eqnarray*}
The coefficient of $x_{ij}^{s+1-m}x_{il}^mx_{kj}^mx_{kl}^{s+1-m}$ is
\begin{eqnarray*}&&(-q^2)^m\begin{pmatrix}s\\m-1\end{pmatrix}_{q^{-4}}q^{4m(m-s-1)}
+(-q^2)^m\begin{pmatrix}s\\m\end{pmatrix}_{q^{-4}}q^{4m(m-s-1)}\\\nonumber
&+&(-q^2)^{m-1}\begin{pmatrix}s\\m-1\end{pmatrix}_{q^{-4}}q^{4m(m-s-1)+2}+(-q^2)^m
\begin{pmatrix}s\\m-1\end{pmatrix}_{q^{-4}}q^{4(m-1)(m-s-1)}\end{eqnarray*}
which can be simplified into
$\begin{pmatrix}s+1\\m\end{pmatrix}_{q^{-4}}q^{4m(m-s-1)}$.

\qed
\begin{Rem}
The third equation in the above lemma can be re-written as
$$\left[b\begin{pmatrix}1&0\\0&1\end{pmatrix}\right]^s
=\sum_{m=0}^s(-q^2)^m\begin{pmatrix}s\\m\end{pmatrix}_{q^{-4}}
x\begin{pmatrix}s-m&m\\m&s-m\end{pmatrix}.$$ \end{Rem}

\section{Kashiwara operators}

\medskip

Let $A_{n-1}=(a_{ij})_{i,j=1}^{n-1}$ be the Cartan matrix of type
$A$. The quantized enveloping algebra $U_q(A_{n-1})$ is the unital
associative algebra generated by  $E_i$, $F_i$, $K_i$, $K_i^{-1}$
($1\le i \le n-1$) with relations
\begin{gather*}
K_i K_j =K_j K_i, \quad K_i K_i^{-1} = 1 = K_i^{-1}K_i \\
E_iF_j - F_jE_i = \delta_{ij}\frac{K_i^2-K_i^{-2}}{q^2-q^{-2}}\\
K_i E_j = q^{2a_{ij}}E_jK_i, \quad K_i F_j = q^{-2a_{ij}}F_jK_i\\
\sum_{s=0}^{1-a_{ij}} (-1)^s \begin{pmatrix}1-a_{ij}\\s\end{pmatrix}_{q^2}
E_i^{1-a_{ij}-s}E_j E_i^s = 0 \quad (i\ne j)\\
\sum_{s=0}^{1-a_{ij}} (-1)^s
\begin{pmatrix}1-a_{ij}\\s\end{pmatrix}_{q^2} F_i^{1-a_{ij}-s}F_j
F_i^s = 0 \quad (i\ne j).
\end{gather*}
As is well known, $U_q(A_{n-1})$ is a Hopf algebra. We take the
following co-multiplication $\Delta$, the co-unit $\epsilon$, and
the antipode $S$ respectively defined by
\begin{align}
\begin{aligned}
& \Delta(K_i) = K_i \otimes K_i, \\
& \Delta(E_i) = E_i \otimes 1 + K_i^2 \otimes E_i,
\quad \Delta(F_i) = F_i \otimes K_i^{-2} + 1 \otimes F_i, \\
& \epsilon(K_i)=1, \quad
\epsilon(E_i) = \epsilon(F_i) = 0, \\
& S(K_i) = K_i^{-1}, \quad
S(E_i) = - K_i^{-2} E_i, \quad
S(F_i) = - F_i K_i^2
\end{aligned}
\end{align}

There are two natural left actions of $U_q(A_{n-1})$ on the
quantized function algebra ${\s O}_q(SL(n))$, which correspond to
left and right translations in the classical setting. These
actions are respectively defined, for all $x\in U_q(A_{n-1}), f\in
{\s O}_q(SL(n))$, by
$$R_x(f)=\sum_{(f)}f_{(1)}<f_{(2)}, x>,$$
$$L_x(f)=\sum_{(f)}<f_{(1)}, S(x)>f_{(2)},$$
where $S$ is the antipode of $U_q(A_{n-1})$. Here we have used
Sweedler's notation $\Delta(f)=\sum f_{(1)}\otimes f_{(2)}$. Note
that $L_x$ and $R_x$ act on ${\s O}_q(SL(n))$ as generalized
derivations. More precisely, for $f, g\in {\s O}_q(SL(n))$,
\begin{eqnarray*}
L_x (f g) &=& \sum_{(x)} L_{x_{(1)}}(f)\cdot L_{x_{(2)}}(g), \\
R_x (f g) &=& \sum_{(x)} R_{x_{(1)}}(f)\cdot R_{x_{(2)}}(g).
\end{eqnarray*} Furthermore, the two left actions commute.

A quantum analogue of the Peter-Weyl theorem states the following:
\begin{Thm} As a left $L_{U_q(A_{n-1})}\otimes R_{U_q(A_{n-1})}$-module,
$${\s O}_q(SL(n))\cong \oplus_{\lambda\in P_+} L^*(\lambda)\otimes L(\lambda).$$
Any $u\otimes v\in L^*(\lambda)\otimes L(\lambda)$ can be viewed as
a linear functional on $U_q(A_{n-1})$ as follows:
$$(u\otimes v)(x)=u(x v), \quad \forall x\in U_q(A_{n-1}).$$
\end{Thm}
The theorem follows from the complete reducibility of finite
dimensional $U_q(A_{n-1})$-modules and the fact that all finite
dimensional $U_q(A_{n-1})$-modules are direct sums of submodules of
tensor powers of the natural module. In fact the theorem holds for
all quantized enveloping algebras of finite dimensional simple Lie
algebras. See \cite{nmy} for a proof in the case of type $A$, and
also see \cite{k} for the general case.

The following observation will be useful for the remainder of the
paper. There exists an algebra anti-involutions of $\omega:
U_q(A_{n-1})\rightarrow U_q(A_{n-1})$ given by
\begin{eqnarray*} \omega(E_i)=F_i, \quad  \omega(F_i)=E_i,
\quad  \omega(K_i)=K_i, \quad  \omega(q)=q. \end{eqnarray*} Then
$\theta=\omega\circ S^{-1}$ is an algebra involution of
$U_q(A_{n-1})$. Using $\theta$ we can construct another left
action of $U_q(A_{n-1})$ on ${\s O}_q(SL(n))$:
\begin{eqnarray} U_q(A_{n-1})\otimes {\s O}_q(SL(n)) \rightarrow  {\s
O}_q(SL(n)),&\quad&  x\otimes f \mapsto L_{\theta(x)}(f).
\label{LeftAction}
\end{eqnarray}

This left action is much simpler than $L$ itself. Let us define
Kashiwara operators for this left action. To this end, we need to
have a suitable basis.
\begin{Prop}There exists an basis of the algebra $O_q(M(n))$
consisting  of the elements of the  form
\begin{eqnarray*}&&q^m x\begin{pmatrix}0&\cdots&0&a_{1r}
&a_{1,r+1}&\cdots&a_{1n}\\a_{21}&\cdots&a_{2, r-1}&
a_{2r}&0&\cdots&0\end{pmatrix}\\ \nonumber && \times \Pi M_{ij}
\Pi_{i\ge 3, j} x_{ij}^{a_{ij}},\end{eqnarray*} where
$a_{ij}\in{\mathbb Z}_+$ for all $i,j$, and the $2\times 2$
quantum minors $M_{i j}$ are of the form $\det_q(\{1, 2\}, \ \{i,
j\})$. The product $\Pi M_{ij}$ of quantum minors is arranged
according to the lexicographic order, namely, $M_{ij}\ge M_{st}$
if $j>t$ or $j=t$ and $i\ge s$. The integer $m$ is uniquely
determined by the entries $a_{ij}$ and the $2\times 2$ minors. The
transition matrix between this new basis and the PBW basis
consisting of the modified monomials is of the form:
$$\begin{pmatrix}1&\cdots &q{\mathbb Z}[q]\\
0&1&\cdots\\\cdots&\cdots&\cdots\\0&\cdots&1\end{pmatrix}$$
\end{Prop}

\pof We will build the basis step by step. First note that the
elements $x_{ij}$ for $i\ge 3$ are annihilated by the left actions
of both $E_1$ and  $F_1$, and so they will contribute nothing
under the Kashiwara operators which will be defined later. Hence,
we may put aside the elements $x_{ij}$ for $i\ge 3$ and only
consider the elements $x_{1l}, x_{2l}$ for $l=1,2,\cdots, n$.

For a PBW basis element $x(A)$, \begin{itemize}
\item if $a_{1,
n-1}a_{2n}=0$, we just keep it unchanged.

\item If $a_{1, n-1}\ge a_{2n}\ge 1$, we can rewrite the modified
monomial $x(A)$ as
\begin{eqnarray*}
x(A)&=&q^{a_{2n}(a_{1n}+a_{2,n-1})}D(A)x_{11}^{a_{11}} \cdots
x_{1, n-2}^{a_{1, n-2}}x_{1,n-1}^{a_{1, n-1}-a_{2n}}
x_{1n}^{a_{1n}}\\
&&\times x_{21}^{a_{21}}\cdots x_{2,n-1}^{a_{2, n-1}}M_{n-1,
n}^{a_{2n}} +\sum_{B<A}c_{BA}x(B). \end{eqnarray*} An elementary
computation shows that $c_{BA}\in q{\mathbb Z}[q]$ with the help
of equation (\ref{power}),

\item If $1\le a_{1, n-1}<a_{2n}$, we can rewrite the modified
monomial $x(A)$ as
\begin{eqnarray*}
x(A)&=&q^{a_{1,n-1}(a_{1n}+a_{2,n-1})}D(A)x_{11}^{a_{11}}\cdots
x_{1, n-2}^{a_{1, n-2}}x_{1n}^{a_{1n}} \\
&&\times x_{21}^{a_{21}}\cdots x_{2,n-1}^{a_{2,
n-1}}x_{2n}^{a_{2n}-a_{1, n-1}}M_{n-1, n}^{a_{1, n-1}}
+\sum_{B<A}c_{BA}x(B).
\end{eqnarray*}
By using equation (\ref{power}), a simple computation shows that
$c_{BA}\in q{\mathbb Z}[q]$.
\end{itemize}
This means that we construct a new basis by replacing the sub-word
$x_{1,n-1}x_{2n}$ by the $2\times 2$ minor $M_{n-1, n}$. Clearly,
the transition matrix between the new basis and the PBW basis of
modified monomials is of the form

$$\begin{pmatrix}1&\cdots &q{\mathbb Z}[q]\\
0&1&\cdots\\\cdots&\cdots&\cdots\\0&\cdots&1\end{pmatrix}.$$

Repeat the above procedure according to the lexicographic order,
replacing the sub-word $x_{1,j}x_{2,k}$ by the $2\times 2$ minor
$M_{jk}$, we get the desired  basis. Note that in each step the
transition between the old basis and the new one is of the form
$$\begin{pmatrix}1&\cdots &q{\mathbb Z}[q]\\
0&1&\cdots\\\cdots&\cdots&\cdots\\0&\cdots&1\end{pmatrix},$$ the
transition matrix between the PBW basis of modified monomials and
our final basis is composition of the matrices of the above form.
This completes the proof. \qed

Denote by $x(A)_1$ the resulting basis element obtained  from
$x(A)$. By doing the same thing to the $i$th and $(i+1)$th rows,
we get a basis $\{x(A)_i| A\in M_n({\mathbb Z}_+)\}$. The order of
$M_n({\mathbb Z}_+)$ gives an order on the basis
$$B_i:=\{x(A)_i|A\in M_n({\mathbb Z}_+)\}.$$
From the transition matrix between $B_i$ and the PBW basis
of modified monomials, we get
$$\overline{x(A)_i}=x(A)_i+\sum_{B<A}e_{BA}x(B)_i$$
with $e_{BA}\in {\mathbb Z}[q, q^{-1}]$.

\begin{Rem}Since the transition matrices of all these various
bases are all of the form
$$\begin{pmatrix}1&\cdots & q{\mathbb Z}[q]\\
0&1&\cdots\\\cdots&\cdots&\cdots \\0&\cdots&1\end{pmatrix},$$ we
have ${\s L}^*=\oplus {\mathbb Z}[q]x(A)=\oplus {\mathbb
Z}[q]x(A)_i$, for all $i$. The canonical basis $B^*$ can be
constructed starting from any basis  $\{x(A)_i|A\in M_n({\mathbb
Z}_+)\}$ for any $i$. More explicitly, for any matrix $A$, the
element $b(A)$ can be written as
$$b(A)=x(A)_i+\sum_{B<A}c_{B,A,i}x(B)_i$$ with $c_{B, A, i}\in
q{\mathbb Z}[q]$. This again justifies the word canonical.
\end{Rem}

Let us consider Kashiwara operators. Note that when examining the
actions of $\tilde{E}_1$ and $\tilde{F}_1$ on $x(A)_1$, we can
ignore the $2\times 2$ minors and those $x_{ij}$ for $i\ge 3$ in
$x(A)_1$. Now the Kashiwara operators $\tilde{E_1}$ and
$\tilde{F_1}$ for the left action (\ref{LeftAction}) are defined
as follows:
\begin{eqnarray*}&&\tilde{E_1}\left(x\begin{pmatrix}0&\cdots&0&a_{1r}
&a_{1,r+1}&\cdots&a_{1n}\\a_{21}&\cdots&a_{2, r-1}&
a_{2 r}&0&\cdots&0\end{pmatrix}\right) \\
&&=\sum_k q^{\sum_{t=1}^{k-1}2a_{2t}}x\begin{pmatrix}0
&\cdots&\quad 1&\cdots& 0& a_{1 r}&a_{1,r+1}&\cdots& a_{1n}
\\a_{21}&\cdots&a_{2k}-1&
\cdots&a_{2, r-1}& a_{2r}&0&\cdots&0\end{pmatrix} \end{eqnarray*}
where the summation is over all the $k$ such that $a_{2k}\ge 1$.
In particular, if $a_{2 r}\ge 1$, the right hand side of the above
equation will contain the term
\[ q^{\sum_{t=1}^{r-1}2a_{2t}}x
\begin{pmatrix}0&\cdots&0&a_{1r}+1
&a_{1,r+1}&\cdots&a_{1n}\\a_{21}&\cdots&a_{2, r-1}& a_{2
r}-1&0&\cdots&0\end{pmatrix}.\] Similarly, we also have
\begin{eqnarray*}
&&\tilde{F_1}\left(x\begin{pmatrix}0&\cdots&0&a_{1r}
&a_{1,r+1}&\cdots&a_{1n}\\a_{21}&\cdots&a_{2, r-1}&
a_{2 r}&0&\cdots&0\end{pmatrix}\right) \\
&&=\sum_{k}q^{\sum_{t>k}2a_{1t}}x
\begin{pmatrix}0&\cdots&0&a_{1r} &a_{1,r+1}&\cdots&a_{1
k}-1&\cdots&a_{1n}\\a_{2 1}&\cdots&a_{2, r-1}& a_{2
r}&0&\cdots&1&\cdots&0\end{pmatrix}
\end{eqnarray*}
where the summation is over the $k$ such that $a_{1k}\ge 1$.

\begin{Prop}\label{action}For any $b(A)\in B^*$ and Kashiwara
operators $\tilde{E_i}$ and $\tilde{F_i}$, $\tilde{E_i}(b(A))=0$
if and only if $\tilde{E_i}(x(A)_i)=0$, and $\tilde{F_i}(b(A))=0$
if and only if $\tilde{F_i}(x(A)_i)=0$.
\end{Prop}

\pof The element $b(A)$ can be written as
\begin{equation} b(A)=x(A)_i+\sum_{B<A}c_{B, A, i}x(B)_i.
\label{b(A)}
\end{equation} We only prove the
statement for $\tilde{E_1}$, the other cases are pretty much the
same. The equation  $\tilde{E_1}x(A)_1=0$ implies that $x(A)_1$ is
of the form:
$$q^sx\begin{pmatrix}0&\cdots&0&a_{1r}&\cdots&a_{1n}\\
0&\cdots &0&0&\cdots&0\end{pmatrix}\times \text{ some }2\times 2
\text{ minors}\times \Pi_{i\ge 3,j}x_{ij}^{a_{ij}}.$$

Note that if we set $i=1$ in equation (\ref{b(A)}), then each
$x(B)_1$ on the right hand side is associated with a $2\times n$
matrix with the second row being zero. This can be seen by
considering $\overline{x}(A)_1$. Using the commutation relations
among the $x_{i j}$ we see that $\overline{x}(A)_1$ can be
expressed as a combination of such $x(B)_1$'s. Hence by the
construction of the dual canonical basis, all of the terms in
$b(A)$ are of the form $x(B)_1$ associated with a $2\times n$
matrix with the second row being zero. Such terms are all
annihilated by $\tilde{E_1}$.

Conversely, if $\tilde{E_1}x(A)_1\ne 0$, then $x(A)_1$ is of the
form
$$q^sx\begin{pmatrix}0&\cdots&0&a_{1r}&\cdots&a_{1n}\\a_{21}&\cdots
&\cdots&a_{2r}&\cdots&0\end{pmatrix}\times \text{ some }2\times 2
\text{ minors}\times \Pi_{i\ge 3,j}x_{ij}^{a_{ij}},$$ with some
$a_{2j}$ nonzero. Without losing generality, we may assume that
$a_{21}\ne 0$. Then the leading term of $\tilde{E_1}x(A)_1$ is of
the form \[q^s x\begin{pmatrix}
1&\cdots&0&a_{1r}&\cdots&a_{1n}\\
a_{21}-1&\cdots &0&a_{2r}&\cdots&0\end{pmatrix}\times \text{ some
} 2\times 2 \text{ minors}\times \Pi_{i\ge 3,j}x_{ij}^{a_{ij}}.\]
In order to obtain a term of this form by applying $\tilde{E_1}$
to some $x(B)_1$, then $x(B)_1$ must be of the form
$$q^sx\begin{pmatrix}1&\cdots&\cdots&a_{1r}&\cdots&a_{1n}\\
a_{21}-1&\cdots&a_{2k}+1&\cdots&\cdots&\cdots
\end{pmatrix}\times \text{some }2\times 2\text{ minors}\times
\Pi_{i\ge 3,j}x_{ij}^{a_{ij}},$$ where $s\in {\mathbb Z}$. As we
have already discussed, such terms can never appear in $b(A)$.
Therefore,  the leading term of $\tilde{E_1}x(A)_1$ can not be
cancelled by terms coming from $\tilde{E_i}x(B)_i$ for $B<A$.
Hence $\tilde{E_i}(b(A)\ne 0$. \qed

The Kashiwara operators are compatible with the action
(\ref{LeftAction}) of Chevalley generators on the quantized
function algebra in the following sense:

\begin{Lem} \label{K-operators} Let $x\in {\mathbb C}_q(SL(n))$.
Then $L_{\theta(E_i)}(x)=0$ if and only if $\tilde{E_i}(x)=0$.
Similarly, $L_{\theta(F_j)}(x)=0$ if and only if
$\tilde{F_j}(x)=0$. \end{Lem}

\pof We shall only prove the case when $i=j=1$. The other cases
are exactly the same.  For convenience we simply write
$L_{\theta(u)}(x)$ as $u(x)$. For a basis element
\begin{eqnarray*}&&x(A)_1=\\
&&q^sx\begin{pmatrix}0&\cdots&0&a_{1r}&\cdots&a_{1n}\\
a_{21}&\cdots &a_{2, r-1}&a_{2r}&\cdots&0\end{pmatrix}\times
\text{some }2\times2 \text{ minors}\times \Pi_{i\ge
3,j}x_{ij}^{a_{ij}}, \end{eqnarray*} both $E_1(x(A)_1)=0$ and
$\tilde{E_1}(x(A)_1)=0$ equivalent to $a_{2i}=0$ for all $i$.

Write $x=\sum_A c_A x(A)_1$. Suppose $E_1(x)=0$. Choose a term
$x(A)_1$ with nonzero coefficient $c_A$ and is maximal with
respect to this property. If $E_1(x(A)_1)\ne 0$, the same argument
as in the above Proposition shows that the leading term in
$E_1(x(A)_1)$ can not be cancelled by the rest of the terms in the
expression of $E_1(x)$ which is a contradiction! Hence, all terms
in the expression of $x$ are killed by $E_1$ which is equivalent
to say that all terms are killed by $\tilde{E_1}$.

The statement for $F_1$ and $\tilde{F_1}$ can be proved similarly.
\qed

From the definition of the Kashiwara operators, we can easily
deduce that

\begin{Prop} For any $b^*\in B^*$, if $\tilde{E_i}(b^*)\ne 0$, then there exists an
element  $b_1^*\in B^*$ such that
$\tilde{E_i}(b^*)=b_1^*\quad mod q{\s L}^*$.
If $\tilde{F_i}(b^*)\ne 0$, then there exists an element $b_2^*\in B^*$ such that
$\tilde{E_i}(b^*)=b_2^*\quad mod q{\s L}^*$.\end{Prop}

\pof If  $\tilde{E_i}(b^*)\ne 0$, then the coefficient of the
leading term in $\tilde{E_i}(b^*)\ne 0$ is 1, which means
there exist a $x(B)_1$ such that
$$\tilde{E_i}(b^*)=x(B)_1\quad mod q{\s L}^*.$$
Hence,  $\tilde{E_i}(b^*)=b(B)\quad mod q{\s L}^*$.

Similarly, we can prove the statement for $\tilde{F_1}$. \qed

\begin{Rem}
Applying the automorphism $\sigma$, we can define the Kashiwara
operators for the right action and the same statements for the
left action hold for the right action.
\end{Rem}

\section{Invariant subalgebras and string bases}\label{inv}

To simplify notations, we set $G=SL(n)$ and ${\mathfrak
g}={\mathfrak{sl}}_n$ in this section. Any subset $S$ of the
generators $\{E_i, \, F_i, \, K_i^{\pm 1} \mid i=1, 2, \dots,
n-1\}$ generate a subalgebra $U_S$ of $U_q({\mathfrak g})$. Note
that $\theta(S)$ also generates a subalgebra
$U_{\theta(S)}=\theta(U_S)$ of $U_q({\mathfrak g})$.
\begin{Def} ${\s O}_q(G)^S:=\{f\in
{\s O}_q(SL(n)) \mid L_x(f)=\epsilon(x) f,  \forall x\in
U_S\}$.\end{Def} This is the subalgebra of ${\s O}_q(G)$ consisting
of the elements which are invariant under the left action $L$ of
$U_S$. As the left translation $L$ and right translation $R$
commute, ${\s O}_q(G)^S$ forms a left $U_q(A_{n-1})$-module under
$R$. The following lemma easily follows from the definition of ${\s
O}_q(G)^S$ and Lemma \ref{K-operators}.
\begin{Lem} If $E_i, F_j\in \theta(S)$, then
$$\tilde{E}_i (f) = 0, \quad \tilde{F}_j (f) = 0,
\quad \forall f\in {\s O}_q(G)^S, $$ where $\tilde{E}_i$ and
$\tilde{F}_j$ are the Kashiwara operators associated with $E_i$
and $F_j$.
\end{Lem}
In an analogous way, we shall also define the subset of dual
canonical basis elements which are invariant with respect to
$U_S$:
\begin{Def} $(\bar{B}^*)^S:=\{b \in \bar{B}^* \mid
L_x(b)=\epsilon(x) b,  \forall x\in U_S\}$.\end{Def}

The following theorem is the main result of this paper.
\begin{Thm} \label{main} The subset $(\bar{B}^*)^S$ of $\bar{B}^*$ forms
a basis of the subalgebra of invariants ${\s O}_q(G)^S$.
\end{Thm}
\pof  Let us express an element $f$ in the subalgebra of
invariants as \[f =\sum_{b\in B^*} c_b b.\] We claim that $b\in
(\bar{B}^*)^S$ if the coefficient $c_b\ne 0$.

Let $e$ be an $E_i$ or $F_j$ belonging to $\theta(S)$. Consider a
term $c_b b$ in $f$ with nonzero coefficient $c_b$, and with $b$
being maximal with respect to the partial order used for defining
the dual canonical basis. Then $L_{\theta(e)} (b)$ must be zero.
Otherwise, the corresponding Kashiwara operator $\tilde{e}$ acts
to produce $\tilde{e}(b)=b' \ mod q{\s L}^*$. It follows from a
similar argument as that in Proposition \ref{action} that $b'$
will not be cancelled by other terms in $\tilde{e}(f)$. \qed

\bigskip

Let $\langle (\bar{B}^*)^S \rangle_+ := {\mathbb Z}_+[q,
q^{-1}](\bar{B}^*)^S$ denote the ${\mathbb Z}_+[q, q^{-1}]$ span
of $(\bar{B}^*)^S$. Then for all $f, g\in \langle (\bar{B}^*)^S
\rangle_+$,
\begin{equation} \label{string-property}
f g \in \langle (\bar{B}^*)^S \rangle_+, \quad \mbox{and} \quad
R_{E_i}(f) \in \langle (\bar{B}^*)^S \rangle_+, \forall i.
\end{equation}
The first property follows from Theorem \ref{positivity}. The
second property can be seen by embedding $B^*$ into a dual
canonical basis of the strictly lower triangular subalgebra
$U_q(A_{2n -1}^-)$ (i.e., the subalgebra generated by the
Chevalley generators $F_i$'s) of $U_q(A_{2n -1})$.  The embedding
was established in \cite{jz}.  It is known \cite{lu2} that the
dual canonical basis of $U_q(A_{2n -1}^-)$ satisfies the second
relation of (\ref{string-property}).

Recall that in \cite{bz1} Berenstein and Zelevinsky introduced the
notion of a string basis, which is a canonical basis with the
${\mathbb Z}_+[q, q^{-1}]$ span satisfying properties analogous to
(\ref{string-property}).
\begin{Lem} $(\bar{B}^*)^S$ forms a string basis for ${\s O}_q(G)^S$.
\end{Lem}

Let us now consider in some detail examples of the subset $S$ and
the associated invariant subalgebra ${\s O}_q(G)^S$.

\subsection{Example: $S=\{F_i\mid i=1, 2,\dots,
n-1\}$}\label{Example1} In this case we denote the subalgebra
$U_S$ by $U_q({\mathfrak n}_{-})$, and ${\s O}_q(G)^S$ by ${\s
O}_q(G)^{U_q({\mathfrak n}_-)}$. Denote by ${\s
O}_q(G)^{U_q({\mathfrak  n}_-)}_\lambda$ the subspace consisting
of elements with $L_{U_q(A_{n-1})}$-weight $-\lambda$ (the
negative of $\lambda$). The result below easily follows from the
quantum Peter-Weyl theorem.
\begin{Thm} \label{BW} The subspace
${\s O}_q(G)^{U_q({\mathfrak  n}_-)}_{\lambda}$, $\lambda\in P_+$,
forms an irreducible left $U_q(A_{n-1})$-submodule with highest
weight $\lambda$ under the action $R$. Furthermore, we have the
$U_q(A_{n-1})$-module isomorphism \[{\s O}_q(G)^{U_q({\mathfrak
n}_-)} = \oplus_{\lambda\in P_+} {\s O}_q(G)^{U_q({\mathfrak
n}_-)}_{\lambda}.\]
\end{Thm}
The first statement of the theorem is a variation of the quantum
analogue \cite{nmy, gz} of the celebrated Borel-Weil theorem.

Let $(\bar{B}^*)_\lambda^{U_q(A_{n-1})}$ denote the subset of
$(\bar{B}^*)^{U_q(A_{n-1})}$ consisting of elements which have
weight $-\lambda$ under the $L$-action of $U_q(A_{n-1})$. Set
\[
\langle(\bar{B}^*)^{U_q({\mathfrak n}_-)}_\nu\rangle_+ = {\mathbb
Z}_+[q, q^{-1}](\bar{B}^*)^{U_q({\mathfrak n}_-)}_\nu, \quad
\nu\in P_+.
\]
Since every element of $\bar{B}^*$ has a definite
weight, we have the following result.
\begin{Thm}\label{irrep} For every $\lambda\in
P_+$, $(\bar{B}^*)^{U_q({\mathfrak n}_-)}_\lambda$ forms a basis
of ${\s O}_q(G)^{U_q({\mathfrak n}_-)}_{\lambda}$. Furthermore,
for all $f\in\langle (\bar{B}^*)^{U_q({\mathfrak
n}_-)}_\lambda\rangle_+$ and $g\in\langle
(\bar{B}^*)^{U_q({\mathfrak n}_-)}_\mu\rangle_+$, we have $f g$
belongs to $\langle(\bar{B}^*)^{U_q({\mathfrak
n}_-)}_{\lambda+\mu}\rangle_+.$
\end{Thm}

By Theorem \ref{BW}, every finite dimensional irreducible
$U_q(A_{n-1})$-module can be realized in terms of  ${\s
O}_q(G)^{U_q({\mathfrak  n}_-)}_{\lambda}$ with the appropriate
highest weight $\lambda\in P_+$.  Thus Theorem \ref{irrep} gives
rise to a basis for each finite dimensional irreducible
$U_q(A_{n-1})$-module.

Let $\Lambda_1,\cdots,\Lambda_{n-1}$ be the fundamental dominant
weights. The natural representation is the simple modules
$L(\Lambda_1)$ which is dual to  $L(\Lambda_{n-1})$ and both are
of dimension $n$.

The quantum minor \[\Delta_s:=det_q(\{n-s+1,\cdots,n\},\
\{1,2,\cdots,s\})\] is of weight $\Lambda_s^*$ (the highest weight
of the dual module of $L(\Lambda_s)$) with respect to the
$L$-action of $U_q(A_{n-1})$, and weight $\Lambda_s$ under the
$R$-action of $U_q(A_{n-1})$. It is easy to see that $\Delta_s$ is
annihilated by all $L_{E_i}$ and $R_{E_i}$ for $i=1,2,\cdots,
n-1$, thus is the highest weight vector of the submodule
$L(\Lambda_s)^*\otimes L(\Lambda_s)$ in ${\s O}_q(SL_n)$.

Moreover, for arbitrary dominant weight
$$\lambda=\ell_1\Lambda_1+\ell_2\Lambda+\cdots+
\ell_{n-1}\Lambda_{n-1},$$ the submodule $L(\lambda)^*\otimes
L(\lambda)$ in ${\s O}_q(SL_n)$ is generated by a monomial
\[\Delta_1^{\ell_1}\Delta_2^{\ell_2}\cdots\Delta_{n-1}^{\ell_{
n-1}}\] which is again a highest weight vector. All such monomials
are dual canonical basis elements by \cite{zhanghc}~ Theorem 5.2.

\begin{Rem}\label{canonical}
The basis $\langle(\bar{B}^*)^{U_q({\mathfrak
n}_-)}_{\lambda+\mu}\rangle_+$ seems to deserve the name of
canonical basis as its elements are bar-invariant ${\mathbb Z}[q]$
combinations of PBW basis elements.
\end{Rem}

\subsection{Example:  $S=\{E_i, \ F_i \mid i\in\Theta\}\cup \{K_j^{\pm 1} \mid
j=1, 2, \dots, n-1\}$} Here $\Theta$ is any set of $\{1, 2, \dots,
n-1\}$.  In this case, we shall denote $U_S$ by $U_q({\mathfrak
k})$, which is isomorphic to $U_q({\mathfrak{gl}}_{k_1})\otimes
U_q({\mathfrak {gl}}_{k_2})\otimes\cdots\otimes U_q({\mathfrak
{gl}}_{k_t})$ for some $t\le n$, where $k_i\ge 1$ and $\sum_{i=1}^t
k_i=n-1$. To be more specific, let $N_s= \sum_{r=1}^{s} k_r$, then
we may assume that that $E_i, F_i$ for $N_{s-1} \le i < N_s$ belong
to the $U_q({\mathfrak{gl}}_{k_s})$ subalgebra of $U_q({\mathfrak
k})$. We shall denote ${\s O}_q(G)^S$ by ${\s O}_q(K\backslash G)$,
which is the algebra of functions on some quantum homogeneous space
\cite{gz} determined by $U_q({\mathfrak k})$.

To describe ${\s O}_q(K\backslash G)$, we first consider the
pre-image  of $(\bar{B}^*)^{U_q(A_{n-1})}$ under the map $T$
(defined by (\ref{T})). Needless to say,
$T^{-1}((\bar{B}^*)^{U_q(A_{n-1})})$ is a subset of $B^*$. Each
element of it spans a $1$-dimensional $U_q({\mathfrak k})$-module
with respect to the restriction of the left action $L$.

Under $L_{U_q({\mathfrak k})}$,  the ${\mathbb Q}(q)$ span of $x_{i
j}$, $i, i =1, 2, \dots, n$, decomposes into a direct sum of
submodules
\begin{eqnarray*}
\sum_{i j=1}^n {\mathbb Q}(q)x_{i j} = \bigoplus_{s=1}^t V^n_s,
\quad V^n_s=\sum_{i=N_{s-1}+1}^{N_s} \sum_{j=1}^n {\mathbb Q}(q)x_{i
j}.
\end{eqnarray*}
Each $V^n_s$ is $n$ copies of the natural module over the
$U_q({\mathfrak {gl}}_{k_s})$ subalgebra in $U_q({\mathfrak k})$,
and with all the other $U_q({\mathfrak {gl}}_{k_i})$ acting
trivially. As a $U_q({\mathfrak {gl}}_{k_s})$-module, the
subalgebra of ${\s O}_q(M(n))$ generated by $V_s^n$ is isomorphic
to the direct sum of the $q$-symmetrized tensor powers of $n$
copies of the natural module. The only possible $1$-dimensional
$U_q({\mathfrak {gl}}_{k_s})$-submodules must come from the
determinant module of $U_q({\mathfrak {gl}}_{k_s})$. Therefore,
inside the subalgebra of ${\s O}_q(M(n))$ generated by $V_s^n$,
every $1$-dimensional $U_q({\mathfrak {gl}}_{k_s})$-submodule is a
${\mathbb Q}(q)$-linear combination of products of quantum
determinants of $k_s\times k_s$ sub-matrices of
\begin{equation}
\left[\begin{matrix}
x_{N_{s-1}+1, 1} & x_{N_{s-1}+1, 2} & \dots & x_{N_{s-1}+1, n}\\
x_{N_{s-1}+2, 1} &  x_{N_{s-1}+2, 2} & \dots & x_{N_{s-1}+2, n}\\
\dots & \dots & \dots& \dots \\
x_{N_s, 1} &  x_{N_s, 2} & \dots & x_{N_s+1, n}
\end{matrix}\right]. \label{submatrix}
\end{equation}
Such quantum determinants are quantum minors of the form
$\det_q(I_s, J(s))$, where $I_s=\{N_{s-1}+1, N_{s-1}+2, \dots,
N_s\}$, and $J(s)$ is any subset of $\{1, 2, \dots, n\}$ of
cardinality $k_s$.

It therefore follows that every element of
$T^{-1}((\bar{B}^*)^{U_q(A_{n-1})})$ is a ${\mathbb Q}[q, q^{-1}]$
combination of products of quantum minors $\det_q(I_s, J(s))$,
$s=1, 2, \dots, t$,  and the quantum determinant $\det_q$. Thus
all elements of $(\bar{B}^*)^{U_q(A_{n-1})}$ can be expressed in
terms of the images of the quantum minors under $T$.

The above discussions have established the following result.
\begin{Lem}
The subalgebra ${\s O}_q(K\backslash G)$ of invariants is
generated by  quantum minors of the form $T(\det_q(I_s, J(s)))$,
$s=1, 2, \dots, t$, where $I_s=\{N_{s-1}+1, N_{s-1}+2, \dots,
N_s\}$, and $J(s)$ is any subset of $\{1, 2, \dots, n\}$ of
cardinality $k_s$.
\end{Lem}

\section*{Acknowledgement} The first author would like to
thank Dr. Zhang Guanglian for helpful discussions.

\bibliographystyle{amsplain}

\end{document}